\documentclass[11pt]{article}

\usepackage{amsmath,amssymb,amsfonts,amsthm}
\usepackage{amscd,amsbsy}
\usepackage{mathtools}
\usepackage[mathscr]{eucal}
\usepackage{booktabs,array}
\usepackage{graphicx,epstopdf,floatrow,subfig}
\usepackage{multirow,makecell}
\usepackage{enumerate}
\usepackage{bm}
\usepackage{xcolor}
\usepackage{threeparttable}
\usepackage{comment}
\usepackage{cite}
\usepackage{url}
\usepackage{hyperref}
\hypersetup{hidelinks}

\allowdisplaybreaks[4]

\newtheorem{theorem}{Theorem}

\newtheorem{corollary}[theorem]{Corollary}
\newtheorem{conjecture}[theorem]{Conjecture}
\newtheorem{lemma}[theorem]{Lemma}

\title{\bfseries
Uniqueness of limit cycle{s} in a class of classical trinomial Li\'enard systems
}

\author{
	Hebai Chen$^{1}$\thanks{Email: chen\_hebai@csu.edu.cn},
Yuhuan Lu$^{1}$\thanks{Email: lu\_yuhuan@csu.edu.cn},
Qian Wang$^{1}$\thanks{Email: wang\_qian@csu.edu.cn},
	Rui Zhang$^{1}$\thanks{Email: zhang\_rui@csu.edu.cn}
	\\[0.5em]
	{\footnotesize
		$^1$School of Mathematics and Statistics, HNP-LAMA,
		Central South University,
	}\\
	{\footnotesize
		Changsha, Hunan 410083, P. R. China
	}\\[0.3em]
}

\date{}

\begin{document}
	
	\maketitle
	\hrule
	
		\begin{abstract}
		This paper studies the number of limit cycles of the polynomial Li\'enard system
		$\dot x=y-F(x)$, $\dot y=-x$, where
		$F(x)=a_1x^{2N}+a_2x^{2m+1}+x^{2n+1}$,
		 with integers \(N,m,n\) satisfying $2\le N\le m<n$, $a_1\ne0$ and $a_2\neq 0$.
		It is proved that this system has exactly one limit cycle and that
		every periodic orbit is hyperbolic and attracting
        for  $a_2<0$; this system has
        no limit cycles for $a_2>0$. The   uniqueness proof compares the widths and relative midpoints of the two scalar
		arcs obtained in potential coordinates. 
		This result gives a positive answer to Conjecture 1.4 of
        Cen et al. (Sci. China Math., 2026).
        
		
	\end{abstract}
	
	\vspace{2em}
	\hrule
	\section{Introduction}
	Among the 23 mathematical problems proposed by Hilbert
	in 1900, the second part of the sixteenth problem concerns
	the maximum number and relative positions of limit cycles
	of planar polynomial differential systems of degree $n$;
	see \cite{Arnold,Hilbert,Smale1991}.
	Since linear systems have no limit cycles, the question
	arises for $n\ge2$. Even for quadratic systems, however,
	it remains unresolved despite extensive research.
	Its difficulty was emphasized by Smale~\cite{Smale},
	who regarded it as the most elusive of Hilbert's problems
	apart from the Riemann hypothesis.
    {For surveys of progress on this problem, see \cite{Ilyashenko,Li2003}.}
	
	A more restricted problem, proposed by
	Smale~\cite{Smale1991}, concerns the same questions
	for classical polynomial Li\'enard systems
	\begin{equation}\label{oL}
		\frac{dx}{dt}=y-F(x),\qquad
		\frac{dy}{dt}=-x,
	\end{equation}
	where $F$ is a real polynomial.
	The problem of determining the maximum number and
	relative positions of limit cycles within this class
	is known as Smale's weakened version of Hilbert's
	sixteenth problem. For surveys of results on limit cycles of polynomial Li\'enard systems, see \cite{CTW,LlibreZhang}.
	
	For $\deg F=n$, Lins, de Melo and Pugh~\cite{LMP}
	conjectured that system~\eqref{oL} has at most
	$[(n-1)/2]$ limit cycles, where $[(n-1)/2]$ denotes
	the largest integer less than or equal to $(n-1)/2$.
	They proved the conjecture for $n\le3$ and also
	determined the global phase portraits for $n=2,3$.
	For $n=4$, Li and Llibre~\cite{LL} confirmed
	the conjecture. Subsequently, Chen and Chen~\cite{ChenChen}
	provided the complete bifurcation diagram and all global
	phase portraits in this case.
	
	However, Dumortier, Panazzolo and Roussarie~\cite{DPR}
	constructed systems of the form~\eqref{oL} with at least
	$[(n-1)/2]+1$ limit cycles for every odd $n\ge7$.
	De Maesschalck and Dumortier~\cite{MD} subsequently
	improved this lower bound to $[(n-1)/2]+2$
	for every $n\ge6$. De Maesschalck and Huzak~\cite{MH}
	further established a lower bound of $n-2$
	for every $n\ge6$.
	These results show that the conjecture is false
	for every $n\ge6$.
	
	For $n=5$, Li and Lu~\cite{LiLu} proved that the
	cyclicity of any nondegenerate slow--fast cycle of
	system~\eqref{oL} is at most two,
    and Chen, Chen and Xiao \cite{CCX} proved {in a preprint} that the
	cyclicity of any  slow--fast cycle of
	system~\eqref{oL} is at most two.
	More recently, OpenAI~\cite{OpenAI} claimed in a preprint
	to have proved that system~\eqref{oL} has at most
	two limit cycles in this case.
    {If the claimed upper bound in \cite{OpenAI} is independently verified, it would establish the Lins--de Melo--Pugh conjecture for $\deg F=5$.}
	
	The above results concern the number of limit cycles in terms of the degree of $F$. An alternative approach is to consider the number of nonzero monomial terms in $F$. 
	The idea is closely related to the fewnomial theory initiated by Khovanskii~\cite{Khovanskii}. 
    In algebraic geometry, Khovanskii~\cite{Khovanskii1978} studied Newton polyhedra and toroidal varieties, while Mikhalkin~\cite{Mikhalkin} studied real algebraic curves, moment maps and amoebas.
    Sprott~\cite{Sprott} investigated chaotic polynomial vector fields in $\mathbb{R}^3$ with a minimal number of monomials. Gasull~\cite{Gasull} proposed a variant of Hilbert's sixteenth problem by replacing the degree of a polynomial system with the number of monomial terms. Subsequently, Álvarez, Gasull and Prohens~\cite{AGP} studied the uniqueness of limit cycles for complex differential equations with two monomials, while Gasull, Li and Torregrosa~\cite{GLT} investigated limit cycles of differential equations with three monomials. Further results on systems with two or three monomial nonlinearities were obtained by Gasull, Gin\'e and Torregrosa~\cite{GGT} and Álvarez et al.~\cite{ACG}. Related studies on generalized Abel equations can be found in~\cite{AT}, and a broader discussion of this approach is given by Gasull and Santana~\cite{GS}.
	
	Recently, Cen, Chen, Tian and Zhao~\cite{CCTZ} studied a variant of Lins--de Melo--Pugh's conjecture by considering system~\eqref{oL} with 
	$$ F(x)=\sum_{i=1}^{n}a_ix^{m_i},$$ 
	where $a_i\neq 0$ and $1\leq m_1<m_2<\cdots<m_n$ are integers. Denoting by $\operatorname{Lie}M(n)$ the maximum number of limit cycles of system~\eqref{oL} when $F$ contains $n$ nonzero monomial terms, they proved that $\operatorname{Lie}M(n)\geq n-1$ and $\operatorname{Lie}M(n)=n-1$ for $n=1,2$. For $n=3$, 
    {they considered system~\eqref{oL} with
\begin{equation}\label{tri}
F(x)=a_1x^{m_1}+a_2x^{m_2}+x^{m_3},
\end{equation}
where $a_1a_2\neq0$ and $1\leq m_1<m_2<m_3$ are integers.}
    They obtained the upper bounds for the number of limit cycles in most cases and determined the corresponding global phase portraits. Moreover, they proposed the following conjectures.
	
	\begin{conjecture}[{\cite{CCTZ}}, Conjecture 1.4]\label{conj.1}
		Suppose that $n=3$, $m_1>2$ is even, $m_2$ and $m_3$ are odd, and $a_2<0$. Then the maximum number of limit cycles of system~\eqref{oL} {with $F$ defined by \eqref{tri}} is one.
	\end{conjecture}
	
	\begin{conjecture}[{\cite{CCTZ}}, Conjecture 1.5]\label{conj.2}
		Suppose that $n=3$, $m_1$ and $m_2$ are odd, $m_3$ is even, and $a_1a_2<0$. Then the maximum number of limit cycles of system~\eqref{oL} {with $F$ defined by \eqref{tri}}
        is one.
	\end{conjecture}
	
	\begin{conjecture}[{\cite{CCTZ}}, Conjecture 1.6]\label{conj.3}
		For every positive integer $n$, $ \operatorname{Lie} M(n)=n-1. $
	\end{conjecture}
	
	In the meantime, the first author, Jin, Li and Zhou~\cite{CJLZ} proved that
   the maximum number of limit cycles of system~\eqref{oL} is at least two
  for all even $m_3\geq 6$  when  $n=3$, {$(m_1,m_2)=(1,3)$, and $a_1a_2<0$,} implying that Conjecture~\ref{conj.2} is incorrect.
   Moreover, they proved that the maximum number of limit cycles of system~\eqref{oL} is at least three
   {when $n=3$, $(m_1,m_2,m_3)=(1,9,10)$,
and $a_1a_2<0$,}
implying that Conjecture~\ref{conj.3} is incorrect.
	
	For Conjecture~\ref{conj.1},  since $m_1$ is even
	and $m_2,m_3$ are odd, {we write $(m_1,m_2,m_3)=(2N, 2m+1, 2n+1)$ and rewrite}
   system~\eqref{oL} as
	\begin{equation}\label{system}
		 \frac{dx}{dt}=y-F(x),\qquad \frac{dy}{dt}=-x,
	\end{equation}
	where
	\[
	F(x)=a_1x^{2N}+a_2x^{2m+1}+x^{2n+1}
	\]
	with {integers $N,m,n$ satisfying $2\leq N\leq m<n$, and $a_1a_2\neq0$}.
	Without loss of generality, we assume that $a_1>0$.
	Indeed, {if $a_1<0$, the transformation $(x,y,t)\mapsto(-x,-y,t)$ replaces $a_1$ by $-a_1>0$ and leaves $a_2$ unchanged. Since time is unchanged, this transformation preserves the number and stability of limit cycles.}
	Our main result is as follows.
	\begin{theorem}\label{thm2}
		System~\eqref{system} has no limit cycles if $a_2>0$,
		and has exactly one limit cycle if $a_2<0$.
		In the latter case, the limit cycle is hyperbolic and attracting. 
	\end{theorem}
	
	The following corollary is immediate.
	
	\begin{corollary}
		Conjecture \ref{conj.1} is correct.
	\end{corollary}
	
	The rest of this paper is organized as follows.
In Section~\ref{sec1}, we give an equivalent characterization
of limit cycles in classical Li\'enard systems. In Section~\ref{sec3},
we prove Theorem~\ref{thm2}.

	\section{An equivalent characterization of limit cycles in system (\ref{system})}
	\label{sec1}
	
	For $x\neq0$, applying the Filippov transformation
	$z=x^2/2$  and using $y$ as the
	independent variable to system \eqref{system} gives
	\begin{equation}\notag
		\frac{dz}{dy}=F_i(z)-y,\qquad z>0,
	\end{equation}
	where $i=1,2$ corresponds to the half-planes $x>0$ and $x<0$,
	respectively, and
	\begin{equation}\label{de}
	\begin{aligned}
		F_1(z)
		&=F\left(\sqrt{2z}\right)
		=a_1(2z)^N+a_2(2z)^{(2m+1)/2}+(2z)^{(2n+1)/2},\\
		F_2(z)
		&=F\left(-\sqrt{2z}\right)
		=a_1(2z)^N-a_2(2z)^{(2m+1)/2}-(2z)^{(2n+1)/2}.
	\end{aligned}
	\end{equation}
	Since $a_2\neq0$, we have $F_1, F_2\in C^m \left([0,+\infty)\right)$ and $m\geq2$. Moreover, 
    {each $F_i$ admits a local $C^2$ extension across $z=0$.}
	Note that the transformation $z=x^2/2$ restricts to a
	diffeomorphism from each of the { half-planes} $x>0$ and $x<0$
	onto $z>0$, with inverse $x=\sqrt{2z}$ and $x=-\sqrt{2z}$,
	respectively.
	For $i=1,2$ and $s>0$, 
	let $z_i(y;s)$ denote the solution of the Cauchy initial value problem
	\begin{equation}\label{3}
		\frac{dz}{dy}=F_i(z)-y,\qquad z(F_i(s);s)=s.
	\end{equation}
	Since $F_i\in C^2$, the Cauchy initial value problem~\eqref{3}
	admits a unique local solution $z_i(y;s)$.
	The following two lemmas describe the endpoints of these
	solution arcs and their dependence on the peak height $s$.
	
	\begin{lemma}\label{lem1}
		For each $i=1,2$, $s>0$, and $0\leq h<s$, the positive
		solution arc $z_i(y;s)$ of \eqref{3} intersects the line
		$z=h$ at exactly two points in the $(y,z)$-plane,
		denoted by $(\alpha_i(h,s),h)$ and $(\beta_i(h,s),h)$,
		where
		\[
		-\infty<\alpha_i(h,s)<F_i(h)<\beta_i(h,s)<+\infty.
		\]
		Moreover, $z_i(y;s)$ is strictly increasing on
		$(\alpha_i(h,s),F_i(s))$ and strictly decreasing on
		$(F_i(s),\beta_i(h,s))$. In particular, it attains its
		unique maximum at $y=F_i(s)$, with $z_i(F_i(s);s)=s$.
	\end{lemma}
	
	\begin{proof}
		Fix $i\in\{1,2\}$, $s>0$ and $0\leq h<s$.
		We first claim that the graph of $z_i(y;s)$ intersects
		the curve $y=F_i(z)$ at exactly one point along its
		nonnegative solution arc, namely $(y,z)=(F_i(s),s)$.
		Indeed, differentiating \eqref{3} gives
		\[
		\frac{d^2z_i}{dy^2}
		=F_i'(z_i)\bigl(F_i(z_i)-y\bigr)-1.
		\]
		At the initial point,
		\[
		z_i(F_i(s);s)=s,\qquad
		\frac{dz_i}{dy}(F_i(s);s)=0,\qquad
		\frac{d^2z_i}{dy^2}(F_i(s);s)=-1.
		\]
		Hence $z_i(y;s)$ has a strict local maximum at $y=F_i(s)$.
		More generally, every intersection with $y=F_i(z)$
		is a critical point satisfying $d^2z_i/dy^2=-1$.
		Thus every critical point is a strict local maximum.
		If there were another intersection, the solution would
		have a local minimum between the two strict maxima,
		which is impossible. This proves the claim.
		
		Consequently, the solution is strictly increasing to the
		left of $F_i(s)$ and strictly decreasing to its right,
		with $h<z_i(y;s)\leq s$ along the portion of the arc
		above $z=h$. Set
		\[
		K_i(h,s)=\max_{h\leq z\leq s}|F_i(z)|.
		\]
		By \eqref{3}, we have
		$-K_i(h,s)-y\leq dz_i/dy\leq K_i(h,s)-y$.
		Suppose that the right branch does not intersect $z=h$
		at any finite value of $y$. Boundedness and the continuation
		theorem then imply that it remains above $z=h$ for all
		$y>F_i(s)$. By monotonicity, it has a horizontal asymptote
		$z=z_+^*\geq h$, that is,
		\[
		\lim_{y\to+\infty}z_i(y;s)=z_+^*.
		\]
		However, continuity of $F_i$ gives
		\[
		\frac{dz_i}{dy}
		=F_i(z_i(y;s))-y\longrightarrow-\infty
		\qquad\text{as }y\to+\infty,
		\]
		which is incompatible with this finite limit.
		Similarly, if the left branch did not intersect $z=h$
		at a finite value of $y$, it would have a horizontal
		asymptote $z=z_-^*\geq h$ as $y\to-\infty$, whereas
		\[
		\frac{dz_i}{dy}
		=F_i(z_i(y;s))-y\longrightarrow+\infty,
		\]
		again a contradiction. Hence both branches intersect
		$z=h$ at finite points, denoted by $\alpha_i(h,s)$
		and $\beta_i(h,s)$.
		
		Finally, since $z_i(y;s)$ is strictly increasing to the
		left of $F_i(s)$ and strictly decreasing to its right,
		equation~\eqref{3} gives
		\[
		F_i(h)-\alpha_i(h,s)
		=\left.\frac{dz_i}{dy}\right|_{y=\alpha_i(h,s)}
		\geq0,
		\qquad
		F_i(h)-\beta_i(h,s)
		=\left.\frac{dz_i}{dy}\right|_{y=\beta_i(h,s)}
		\leq0.
		\]
		Equality at either endpoint would give another intersection
		with the curve $y=F_i(z)$ at $(y,z)=(F_i(h),h)$.
		Since $h<s$, this contradicts the first claim. Hence
		\[
		-\infty<\alpha_i(h,s)<F_i(h)<\beta_i(h,s)<+\infty,
		\]
		which	completes the proof.
	\end{proof}
	By Lemma~\ref{lem1}, the arc $z_i(y;s)$ cuts out a segment
	on the line $z=h$. Define its half-length by
	\begin{equation}\label{r}
		r_i(h,s)=\frac{\beta_i(h,s)-\alpha_i(h,s)}{2},
	\end{equation}
	and its relative midpoint with respect to $F_i(h)$ by
	\begin{equation}\label{M}
		\mathcal M_{F_i}(h,s)
		=\frac{\beta_i(h,s)+\alpha_i(h,s)}{2}-F_i(h).
	\end{equation}
	
	\begin{lemma}\label{lem2}
		For each $i=1,2$ and each fixed $h\geq0$, the map
		$s\mapsto r_i(h,s)$ is a strictly increasing $C^1$
		bijection from $(h,\infty)$ onto $(0,\infty)$.
		Consequently, using $r=r_i(h,s)$ in place of $s$,
		the relative midpoint can be regarded as a $C^1$ function
		$\mathcal M_{F_i}(h,r)$ on $[0,\infty)\times(0,\infty)$,
		satisfying
		\begin{equation}\label{8}
			|\mathcal M_{F_i}(h,r)|<r,\qquad
			\left|\frac{\partial}{\partial r}(\mathcal M_{F_i}(h,r))\right|<1,
			\qquad h\geq0,\quad r>0.
		\end{equation}
	\end{lemma}
	\begin{proof}
		Fix $i\in\{1,2\}${. W}e divide the proof into
    the following steps.
		
		\noindent\textbf{Step 1. Derivatives of the endpoints
			$\alpha_i(h,s)$, $\beta_i(h,s)$ with respect to
			the peak height $s$.}
		
		Since $F_i\in C^2([0,\infty))$, the solution  $z_i(y;s)$ of the Cauchy initial value problem
		depends continuously differentiably on the initial  value
		$(F_i(s),s)$, and hence on $s$. Set
		\[
		\xi_i(y;s)=\frac{\partial }{\partial s}(z_i(y;s)).
		\]
		By equation \eqref{3},
		\[
		z_i(y;s)
		=s+\int_{F_i(s)}^y
		\bigl(F_i(z_i(\eta;s))-\eta\bigr)\,d\eta.
		\]
		Differentiating with respect to $s$ gives
		\[
		\begin{aligned}
			\xi_i(y;s)
			={}1-F_i'(s)
			\bigl(F_i(z_i(F_i(s);s))-F_i(s)\bigr)+\int_{F_i(s)}^y
			F_i'(z_i(\eta;s))\xi_i(\eta;s)\,d\eta.
		\end{aligned}
		\]
		Since $z_i(F_i(s);s)=s$,  we  obtain that 
		\[
		\xi_i(y;s)
		=1+\int_{F_i(s)}^y
		F_i'(z_i(\eta;s))\xi_i(\eta;s)\,d\eta,
		\]
		and therefore
		\[
		\frac{d\xi_i}{dy}
		=F_i'(z_i(y;s))\xi_i,
		\qquad
		\xi_i(F_i(s);s)=1.
		\]
		It follows that
		\begin{equation}\label{xi}
		\xi_i(y;s)
		=\exp\left(
		\int_{F_i(s)}^y F_i'(z_i(\eta;s))\,d\eta
		\right)>0.
		\end{equation}
		
		Next, by
		Lemma~\ref{lem1} and equation \eqref{3},
		\[
		(\partial_y z_i)(\alpha_i(h,s);s)
		=F_i(h)-\alpha_i(h,s)>0,
		\qquad
		(\partial_y z_i)(\beta_i(h,s);s)
		=F_i(h)-\beta_i(h,s)<0.
		\]
		Applying the implicit function theorem to
		$z_i(y;s)-h=0$
		at the two endpoints  shows that $\alpha_i(h,s)$
		and $\beta_i(h,s)$ are jointly $C^1$ in $(h,s)$
		for $0\leq h<s$.
		We  therefore differentiate the identities
		\[
		z_i(\alpha_i(h,s);s)=h,\qquad
		z_i(\beta_i(h,s);s)=h
		\]
		with respect to $s$, keeping $h$ fixed, which gives
		\[
		(\partial_y z_i)(\alpha_i(h,s);s)\cdot 
		(\partial_s \alpha_i)(h,s)
		+\xi_i(\alpha_i(h,s);s)=0
		\]
		and 
		\[
		(\partial_y z_i)(\beta_i(h,s);s)\cdot 
		(\partial_s \beta_i)(h,s)
		+\xi_i(\beta_i(h,s);s)=0.
		\]
		Using \eqref{3} and $\xi_i>0$, we obtain
		\begin{equation}\label{4}
			\begin{aligned}
				\frac{\partial}{\partial s}(\alpha_i(h,s))
				=-\frac{\xi_i(\alpha_i(h,s);s)}
				{F_i(h)-\alpha_i(h,s)}<0,\qquad
				\frac{\partial}{\partial s}(\beta_i(h,s))
				=\frac{\xi_i(\beta_i(h,s);s)}
				{\beta_i(h,s)-F_i(h)}>0.
			\end{aligned}
		\end{equation}
		
		\medskip
		\noindent\textbf{Step 2. Strict monotonicity and bijectivity of
			$s\mapsto r_i(h,s)$.}
		
		While keeping $h$ fixed, it follows from \eqref{4} that
		\[
		\frac{\partial }{\partial s}(r_i(h,s))
		=\frac12\left(
		\frac{\partial}{\partial s}(\beta_i(h,s))
		-\frac{\partial}{\partial s}(\alpha_i(h,s))
		\right)>0.
		\]
		Thus $s\mapsto r_i(h,s)$ is $C^1$ and strictly increasing.
		
        By Lemma~\ref{lem1}, for each level $z\in[h,s)$, the two intersections  lie in $[\alpha_i(h,s),\beta_i(h,s)]$ and bracket $F_i(z)$.
        Therefore,  we have that 
		\[
		\alpha_i(h,s)<\min_{z\in[h,s]}F_i(z)
		\leq F_i(s)\leq\max_{z\in [h,s]}F_i(z)<\beta_i(h,s).
		\]
		Integrating \eqref{3} along $(\alpha_i(h,s), F_i(s))$ gives
		\begin{equation}\label{5}
			\begin{aligned}
				s-h
				&=\int_{\alpha_i(h,s)}^{F_i(s)}
				\bigl(F_i(z_i(y;s))-y\bigr)\,dy\\
				&\geq \int_{\alpha_i(h,s)}^{\min_{z\in[h,s]}F_i(z)}
				\left(\min_{z\in[h,s]}F_i(z)-y\right)\,dy
				=\frac12\left(\min_{z\in[h,s]}F_i-\alpha_i(h,s)\right)^2.
			\end{aligned}
		\end{equation}
		Similarly, integrating along $(F_i(s),\beta_i(h,s))$ gives
		\begin{equation}\label{6}
			s-h
			\geq \int_{\max_{z\in[h,s]}F_i(z)}^{\beta_i(h,s)}
			\left(y-\max_{z\in[h,s]}F_i(z)\right)\,dy
			=\frac12\left(\beta_i(h,s)-\max_{z\in[h,s]}F_i(z)\right)^2.
		\end{equation}
		Combining 
		\eqref{5} and \eqref{6} yields
		\[
		0<r_i(h,s)\leq \sqrt{2(s-h)}
		+\frac12\bigl(\max_{z\in [h,s]}F_i(z)-\min_{z\in[h,s]}F_i(z)\bigr).
		\]
		By continuity of $F_i$, we obtain  $\lim_{s\to h^+}r_i(h,s)=0$.
		
		Since $z_i(y;s)\in[h,s]$ for
		$y\in[\alpha_i(h,s),F_i(s)]$, we obtain
		\[
		\begin{aligned}
			s-h
			&\leq \int_{\alpha_i(h,s)}^{F_i(s)}
			\left(\max_{z\in[h,s]}F_i(z)-y\right)\,dy
			\leq \int_{\alpha_i(h,s)}^{\beta_i(h,s)}
			\bigl(\beta_i(h,s)-y\bigr)\,dy\\
			&=\frac12\bigl(\beta_i(h,s)-\alpha_i(h,s)\bigr)^2
			=2r_i^2(h,s).
		\end{aligned}
		\]
		Hence $r_i(h,s)\geq\sqrt{(s-h)/2}$, and therefore
		$\lim_{s\to\infty}r_i(h,s)=\infty$.
		Together with continuity, strict monotonicity, and
		$\lim_{s\to h^+}r_i(h,s)=0$, this proves that
		$s\mapsto r_i(h,s)$ is a bijection from $(h,+\infty)$ onto
		$(0,+\infty)$.
		
		\medskip
		\noindent\textbf{Step 3. Regularity and bounds for the relative midpoint.}
		
		By Steps 1-2, $r_i(h,s)$ is $C^1$ and
		$\partial r_i/\partial s>0$ with $h$ fixed. The implicit function theorem
		therefore gives a $C^1$ inverse $s=S_i(h,r)$ on
		$[0,\infty)\times(0,\infty)$. At $h=0$, we apply the theorem
		using a local $C^2$ extension of $F_i$ across $z=0$.
		Consequently,
		\begin{equation}\notag
			\mathcal M_{F_i}(h,r)
			=\left.\frac{\alpha_i(h,s)+\beta_i(h,s)}{2}
			\right|_{s=S_i(h,r)}-F_i(h)
		\end{equation}
		is  $C^1$.
		
		Since $\alpha_i(h,s)<F_i(h)<\beta_i(h,s)$, it follows from \eqref{r} and \eqref{M} that $|\mathcal M_{F_i}(h,r)|<r$.
		Differentiating with respect to $r$ at fixed $h$ and using
		$r_i(h,S_i(h,r))=r$, we obtain
		\begin{equation}\label{22}
		\begin{aligned}
			\frac{\partial}{\partial r}(\mathcal M_{F_i}(h,r))
			&=\left.
			\frac{
				\dfrac{\partial}{\partial s}
				\left(\dfrac{\alpha_i(h,s)+\beta_i(h,s)}{2}-F_i(h)\right)
			}{
				\dfrac{\partial r_i}{\partial s}(h,s)
			}
			\right|_{s=S_i(h,r)}\\
			&=\left.
			\frac{
				\dfrac{\partial}{\partial s}(\beta_i(h,s))
				+\dfrac{\partial}{\partial s}(\alpha_i(h,s))
			}{
				\dfrac{\partial}{\partial s}(\beta_i(h,s))
				-\dfrac{\partial}{\partial s}(\alpha_i(h,s))
			}
			\right|_{s=S_i(h,r)}.
		\end{aligned}
		\end{equation}
        {By \eqref{4}, we have
\[
    \frac{\partial}{\partial s}(\alpha_i(h,s))<0
    <\frac{\partial}{\partial s}(\beta_i(h,s)),
\]
and hence
\[
    \left|
        \frac{\partial}{\partial s}(\beta_i(h,s))
        +\frac{\partial}{\partial s}(\alpha_i(h,s))
    \right|
    <
    \frac{\partial}{\partial s}(\beta_i(h,s))
    -\frac{\partial}{\partial s}(\alpha_i(h,s)).
\]
Together with \eqref{22}, this gives
\[
    \left|
        \frac{\partial}{\partial r}(M_{F_i}(h,r))
    \right|<1.
\]
This proves \eqref{8}.}
	\end{proof}
	Since $F_1(0)=F_2(0)=0$, the relative midpoints at $h=0$
	are the actual midpoints of the intervals cut out on the
	$y$-axis. Define
	\begin{equation}\notag
		\Delta(r)=\mathcal M_{F_1}(0,r)-\mathcal M_{F_2}(0,r),
		\qquad r>0.
	\end{equation}
	By Lemma~\ref{lem2},
	$\Delta\in C^1((0,\infty))$.
	\begin{theorem}\label{thm1}
		The limit cycles of system~\eqref{system} are in one-to-one
		correspondence with the isolated positive zeros of $\Delta$.
		The limit cycle corresponding to an isolated positive zero
		$r_*$ is hyperbolic if and only if $\Delta'(r_*)\neq0$.
		It is asymptotically stable if $\Delta'(r_*)>0$ and
		unstable if $\Delta'(r_*)<0$.
	\end{theorem}
	
	\begin{proof}
		We first show that every positive zero of $\Delta$ gives a
		periodic orbit. Let $r_*>0$ satisfy $\Delta(r_*)=0$.
		By Lemma~\ref{lem2}, for each $i=1,2$, the half-length $r_*$
		uniquely determines the peak height $S_i(0,r_*)$ and hence
		the endpoint ordinates $\alpha_i(0,S_i(0,r_*))<0$ and
		$\beta_i(0,S_i(0,r_*))>0$ by Lemma \ref{lem1}.
		
		Since $F_i(0)=0$, it follows from \eqref{r} and \eqref{M} that
		\[
		\begin{aligned}
			\alpha_i\bigl(0,S_i(0,r_*)\bigr)
			=\mathcal M_{F_i}(0,r_*)-r_*,\qquad
			\beta_i\bigl(0,S_i(0,r_*)\bigr)
			=\mathcal M_{F_i}(0,r_*)+r_*.
		\end{aligned}
		\]
		Since $\Delta(r_*)=0$, we have
		$\mathcal M_{F_1}(0,r_*)=\mathcal M_{F_2}(0,r_*)$. Hence
		\[
		\begin{aligned}
			\alpha_1\bigl(0,S_1(0,r_*)\bigr)
			&=\alpha_2\bigl(0,S_2(0,r_*)\bigr)=:\alpha,\\
			\beta_1\bigl(0,S_1(0,r_*)\bigr)
			&=\beta_2\bigl(0,S_2(0,r_*)\bigr)=:\beta.
		\end{aligned}
		\]
		Thus, under the inverse transformations $x=\sqrt{2z}$ and
		$x=-\sqrt{2z}$, respectively, the corresponding orbit arcs in the two half-planes
		have the same endpoints $(0,\alpha)$ and $(0,\beta)$. Since $\dot y=-x$, the right arc starts from $(0,\beta)$ and arrives at
		$(0,\alpha)$, and the left arc runs from $(0,\alpha)$ to
		$(0,\beta)$. 
        {Set $s_i=S_i(0,r_*)$ for $i=1,2$. By \eqref{3} and $F_i(0)=0$, we have
\[
    z_i'(\alpha;s_i)=-\alpha>0,\qquad
    z_i'(\beta;s_i)=-\beta<0.
\]
Consequently,
\[
    z_i(y;s_i)=(-\alpha)(y-\alpha)+o(y-\alpha)
    \quad\text{as }y\to\alpha^+,
\]
and
\[
    z_i(y;s_i)=\beta(\beta-y)+o(\beta-y)
    \quad\text{as }y\to\beta^-.
\]
Since $|dt|=|dy|/\sqrt{2z_i(y;s_i)}$, these estimates
show that both endpoints are reached in finite time.
Moreover, $\dot{x}=y\neq0$ at $(0,\alpha)$ and $(0,\beta)$,
so the trajectories cross the $y$-axis transversely.
By the uniqueness of solutions, the two arcs join to form
a periodic orbit of system~\eqref{system}, which meets
the positive and negative $y$-axes at $(0,\beta)$
and $(0,\alpha)$, respectively.}
		
		Suppose further that $r_*$ is an isolated zero of $\Delta$.
        Then $\Delta(r)\neq0$ for all \(r\) sufficiently close to \(r_*\) with \(r\neq r_*\).
		 For such $r$, the endpoint formulas give
		\[
		\begin{aligned}
			\alpha_1\bigl(0,S_1(0,r)\bigr)-\alpha_2\bigl(0,S_2(0,r)\bigr)
			&=\Delta(r)\neq0,\\
			\beta_1\bigl(0,S_1(0,r)\bigr)-\beta_2\bigl(0,S_2(0,r)\bigr)
			&=\Delta(r)\neq0.
		\end{aligned}
		\]
		Thus the two arcs cannot form a periodic orbit for any such $r$. Hence the periodic
		orbit corresponding to $r_*$ is isolated, and thus is a limit cycle.
		
		Conversely, let $\Gamma$ be a limit cycle of system
		\eqref{system}. Since $\Gamma$ surrounds the origin, it is easy to see that $\Gamma$ meets each of the positive
		and negative $y$-axes exactly once, at $(0,\beta)$ and
		$(0,\alpha)$, respectively, with $\alpha<0<\beta$.
		Under the transformation $z=x^2/2$, the portions of $\Gamma$
		in the two half-planes give scalar solution arcs $z_1$ and $z_2$ with the same
		endpoint ordinates $\alpha$ and $\beta$ at $z=0$. Thus they
		have the same half-length and midpoint:
		\[
		r_*=\frac{\beta-\alpha}{2},
		\qquad
		\mathcal M_{F_1}(0,r_*)=\mathcal M_{F_2}(0,r_*)
		=\frac{\alpha+\beta}{2}.
		\]
		Consequently, $\Delta(r_*)=0$.  For $i=1,2$, define $p_i$ by mapping the positive
		ordinate $\beta_i(0,s)$ to the negative ordinate
		$\alpha_i(0,s)$ of the same solution arc $z_i(y;s)$
		at $z=0$. That is,
		\[
		p_i\bigl(\beta_i(0,s)\bigr)=\alpha_i(0,s),\qquad s>0.
		\] Thus,
		for $r$ near $r_*$,
		\begin{equation}\label{9}
			p_i\bigl(r+\mathcal M_{F_i}(0,r)\bigr)
			=\mathcal M_{F_i}(0,r)-r.
		\end{equation}
        {By \eqref{8},
\[
    \frac{d}{dr}\bigl(r+M_{F_i}(0,r)\bigr)
    =1+\frac{\partial}{\partial r}M_{F_i}(0,r)>0,
\]
whereas
\[
    \frac{d}{dr}\bigl(M_{F_i}(0,r)-r\bigr)
    =\frac{\partial}{\partial r}M_{F_i}(0,r)-1<0.
\]
Thus $r\mapsto r+M_{F_i}(0,r)$ is strictly increasing,
and \eqref{9} shows that each $p_i$ is a $C^1$
diffeomorphism onto its image with negative derivative.}
        Following the flow first through $x>0$ and
		then through $x<0$, we obtain the local Poincar\'e return map $P=p_2^{-1}\circ p_1$ of system \eqref{system}
		near $(0,\beta)$.

		Since $\Gamma$ is a limit cycle, $\beta$ is an isolated
		fixed point of $P$. Suppose that $r_*$ is not an isolated
		zero of $\Delta$. Then there is a sequence of zeros
		$r_n\neq r_*$ of $\Delta$ with $r_n\to r_*$. For all
		sufficiently large $n$, the ordinate
		$r_n+\mathcal M_{F_1}(0,r_n)$ is a fixed point of $P$.
		By \eqref{8}, the map $r\mapsto r+\mathcal M_{F_1}(0,r)$
		is continuous and strictly increasing. Thus these fixed
		points are distinct from $\beta$ and converge to
		$\beta=r_*+\mathcal M_{F_1}(0,r_*)$, contradicting the
		isolation of $\beta$. Hence $r_*$ is an isolated zero
		of $\Delta$.

		Finally, let $r_*$ be an isolated positive zero of $\Delta$
		and let $\beta_*$ be the common positive endpoint ordinate.
		Differentiating \eqref{9}  at $r=r_*$ gives
		\[
		p_i'(\beta_*)=\frac{\frac{\partial }{\partial r}(\mathcal M_{F_i}(0,r_*))-1}{\frac{\partial }{\partial r}(\mathcal M_{F_i}(0,r_*))+1}.
		\]
		Consequently,
		\[
		P'(\beta_*)-1
		=-\frac{2\Delta'(r_*)}{\left(\frac{\partial }{\partial r}(\mathcal M_{F_1}(0,r_*))+1\right)\left(1-\frac{\partial }{\partial r}(\mathcal M_{F_2}(0,r_*))\right)}.
		\]
        {By \eqref{8},
\[
    -1<
    \frac{\partial}{\partial r}M_{F_i}(0,r_*)
    <1.
\]
Hence both factors in the denominator of the preceding
expression are positive. Moreover,
\[
    P'(\beta_*)
    =\frac{p_1'(\beta_*)}{p_2'(\beta_*)}>0.
\]
Consequently, $P'(\beta_*)\neq1$ if and only if
$\Delta'(r_*)\neq0$. If $\Delta'(r_*)>0$, then
$0<P'(\beta_*)<1$, and the corresponding limit cycle
is hyperbolic and asymptotically stable.
If $\Delta'(r_*)<0$, then $P'(\beta_*)>1$,
and the corresponding limit cycle is hyperbolic and unstable.}
	\end{proof}

\section{Proof of Theorem~\ref{thm2}}\label{sec3}
{Under the identification $(m_1,m_2,m_3)=(2N,2m+1,2n+1)$,
\cite[Theorem~1.3(iii.1)]{CCTZ} implies that}
system~\eqref{system}
has no limit cycles if $a_2>0$, whereas it has at least one limit cycle if  $a_2<0$.
We therefore restrict our attention to the case \(a_2<0\). It remains to prove that
\(\Delta'(r)>0\) whenever  \(r>0\)
and \(\Delta(r)=0.\)
Since \(\Delta\) is of class \(C^1\), this property implies that \(\Delta\) has at most one positive zero. Consequently, system~\eqref{system} has a unique limit cycle, which is hyperbolic and asymptotically stable by Theorem~\ref{thm1}.

In equations~\eqref{de}, write
\[
F_1(z)=E(z)+H(z),\qquad F_2(z)=E(z)-H(z),
\]
where
\[
E(z)=a_1(2z)^N,\qquad
H(z)=a_2(2z)^{(2m+1)/2}+(2z)^{(2n+1)/2}.
\]
In what follows, we have $a_1>0$ and $a_2<0$.

\begin{lemma}\label{lem9}
Consider  $F_1$ and $F_2$ as defined in \eqref{de}, and fix
	$h_1,h_2\geq0$.  The following results hold:
	\begin{enumerate}
		\item[(i)] If 
		\(F_1(h_1+z)-F_1(h_1)>F_2(h_2+z)-F_2(h_2)\)
		for all \(z>0,\)
		then
		\(\mathcal M_{F_1}(h_1,r)>\mathcal M_{F_2}(h_2,r)\)
		for all \(r>0.\)
		
		\item[(ii)] If
		\(F_1(h_1+z)-F_1(h_1)>F_2(h_2)-F_2(h_2+z)\)
	for all \(z>0,\)
		then
		\(\mathcal M_{F_1}(h_1,r)+\mathcal M_{F_2}(h_2,r)>0\)
		for all \(r>0.\)
	\end{enumerate}
\end{lemma}
\begin{proof}
Fix $r>0$. By Lemma~\ref{lem2}, for each $i=1,2$ and fixed
$h_i\geq0$, there is a unique peak height $s_i=S_i(h_i,r)>h_i$
such that $r_i(h_i,s_i)=r$. Recall that $z_i(y;s_i)$ is the corresponding
solution of the Cauchy initial value problem~\eqref{3}, that
$(\alpha_i(h_i,s_i),h_i)$ and $(\beta_i(h_i,s_i),h_i)$ are its
endpoints on $z=h_i$, and that $\mathcal M_{F_i}(h_i,r)$ is
the  relative midpoint of these endpoints relative to $F_i(h_i)$.

Set
\[
Z_i(y)=z_i\bigl(y+F_i(h_i);s_i\bigr)-h_i.
\]
Then
\[
Z_i'(y)=F_i\bigl(h_i+Z_i(y)\bigr)-F_i(h_i)-y.
\]
The translated arc $Z_i(y)$ intersects the line $z=0$ exactly at
$(\bar\alpha_i,0)$ and $(\bar\beta_i,0)$, where
\(\bar\alpha_i=\alpha_i(h_i,s_i)-F_i(h_i)\) and
\(\bar\beta_i=\beta_i(h_i,s_i)-F_i(h_i).\)
By \eqref{r} and \eqref{M}, together with the endpoint
inequalities in Lemma~\ref{lem1}, we have
\[
\bar\alpha_i=\mathcal M_{F_i}(h_i,r)-r<0,
\qquad
\bar\beta_i=\mathcal M_{F_i}(h_i,r)+r>0.
\]

Suppose that
$\mathcal M_{F_1}(h_1,r)\leq\mathcal M_{F_2}(h_2,r)$.
Then $\bar\alpha_1\leq\bar\alpha_2<0<\bar\beta_1\leq\bar\beta_2$,
and hence
\[
Z_1(\bar\alpha_2)\geq Z_1(\bar\alpha_1)
=Z_2(\bar\alpha_2)=0,
\qquad
Z_1(\bar\beta_1)=Z_2(\bar\beta_2)=0\leq Z_2(\bar\beta_1).
\]
For $y\in[\bar\alpha_2,\bar\beta_1]$, set
$w(y)=Z_1(y)-Z_2(y)$. Then
\begin{equation}\label{w}
w'(y)
=F_1(h_1+Z_1(y))-F_1(h_1)
-\bigl(F_2(h_2+Z_2(y))-F_2(h_2)\bigr).
\end{equation}
Since $F_1\in C^2$, by  the mean value theorem, whenever
$w(y)\neq0$, there exists a point $\xi(y)$ between $h_1+Z_1(y)$ and
$h_1+Z_2(y)$ such that
\begin{equation}\label{f}
F_1(h_1+Z_1(y))-F_1(h_1+Z_2(y))
=F_1'(\xi(y))w(y).
\end{equation}
When $w(y)=0$, set $\xi(y)=h_1+Z_2(y)$. 
Moreover, $F_1'(\xi(y))$ is continuous in $y$.

Since $Z_2(y)>0$ for $y\in(\bar\alpha_2,\bar\beta_1)$, from \eqref{w} and \eqref{f},
the hypothesis of (i) gives
\[
\begin{aligned}
	w'(y)-F_1'(\xi(y))w(y)
=F_1(h_1+Z_2(y))-F_1(h_1)-\bigl(F_2(h_2+Z_2(y))-F_2(h_2)\bigr)>0.
\end{aligned}
\]
Therefore, on $(\bar\alpha_2,\bar\beta_1)$,
\[
\begin{aligned}
	\frac{d}{dy}\left(
	\exp\left(-\int_{\bar\alpha_2}^{y}
	F_1'(\xi(\eta))\,d\eta\right)w(y)
	\right)
	=\exp\left(-\int_{\bar\alpha_2}^{y}
	F_1'(\xi(\eta))\,d\eta\right)
	\bigl(w'(y)-F_1'(\xi(y))w(y)\bigr)
	>0.
\end{aligned}
\]
Integrating over $[\bar\alpha_2,\bar\beta_1]$, we obtain
\[
\exp\left(-\int_{\bar\alpha_2}^{\bar\beta_1}
F_1'(\xi(\eta))\,d\eta\right)w(\bar\beta_1)
-w(\bar\alpha_2)>0.
\]
This contradicts $w(\bar\alpha_2)\geq0$ and
$w(\bar\beta_1)\leq0$, and proves (i).

	For (ii),   if $z(y)$ satisfies $dz/dy=F_2(z)-y$,
	then $z(-y)$ satisfies $dz/dy=-F_2(z)-y$. This reflection
	preserves the half-length and reverses the sign of the
	relative midpoint, so
	\[
	\mathcal M_{-F_2}(h_2,r)=-\mathcal M_{F_2}(h_2,r).
	\]
	The hypothesis of (ii) is precisely the increment inequality
	in (i) with $F_2$ replaced by $-F_2$. The same comparison
	argument therefore gives
	\[
	\mathcal M_{F_1}(h_1,r)>
	\mathcal M_{-F_2}(h_2,r)=-\mathcal M_{F_2}(h_2,r),
	\]
	which proves (ii).
\end{proof}

\begin{lemma}\label{lem10}
		Consider  $F_1$ and $F_2$,  defined by \eqref{de}.
		If $0<u<v$ and $F_1'(u)=F_2'(v)$, then
		$F_1''(u)>F_2''(v)$.
	\end{lemma}
\begin{proof}
	Since $F_1=E+H$ and $F_2=E-H$,  it follows from
	$F_1'(u)=F_2'(v)$  that
	\begin{equation}\label{a}
		\begin{aligned}
		a_2(2m+1)2^{m-\frac12}
			\bigl(u^{m-\frac12}+v^{m-\frac12}\bigr)
		=E'(v)-E'(u)
			-(2n+1)2^{n-\frac12}
			\bigl(u^{n-\frac12}+v^{n-\frac12}\bigr).
		\end{aligned}
	\end{equation}
	On the other hand,
	\[
	\begin{aligned}
		F_1''(u)-F_2''(v)
		={}&E''(u)-E''(v)\\
		&+a_2(2m+1)(2m-1)2^{m-\frac32}
		\bigl(u^{m-\frac32}+v^{m-\frac32}\bigr)\\
		&+(2n+1)(2n-1)2^{n-\frac32}
		\bigl(u^{n-\frac32}+v^{n-\frac32}\bigr).
	\end{aligned}
	\]
By	using \eqref{a}, we obtain
	\[
	\begin{aligned}
		F_1''(u)-F_2''(v)
		={}&E''(u)-E''(v)+\bigl(E'(v)-E'(u)\bigr)
		\left(m-\tfrac12\right)
		\frac{u^{m-\frac32}+v^{m-\frac32}}
		{u^{m-\frac12}+v^{m-\frac12}}\\
		&+(2n+1)2^{n-\frac12}
		\bigl(u^{n-\frac12}+v^{n-\frac12}\bigr)\\
		&\quad\times\left(
		\left(n-\tfrac12\right)
		\frac{u^{n-\frac32}+v^{n-\frac32}}
		{u^{n-\frac12}+v^{n-\frac12}}
		-
		\left(m-\tfrac12\right)
		\frac{u^{m-\frac32}+v^{m-\frac32}}
		{u^{m-\frac12}+v^{m-\frac12}}
		\right).
	\end{aligned}
	\]

Set $a=u/v\in(0,1)$ and define
\[
g(t)=\frac{t(a^{t-1}+1)}{a^t+1}.
\]
We claim that $g(t)$ is strictly increasing on $[3/2,+\infty)$.
By calculation,
\[
g'(t)=\frac{a^{t-1}}{(a^t+1)^2}G(t),
\qquad
G(t)=1+a+a^t+a^{1-t}+t(1-a)\ln a.
\]
For $t\geq 3/2$, we have $a^{1-t}>1$ and $a^t<a$. Thus
\[
G'(t)=(-\ln a)\bigl(a^{1-t}-1+a-a^t\bigr)>0.
\]
Therefore, to prove that $g$ is strictly increasing on
$[3/2,+\infty)$, it suffices to show that $G(3/2)>0$.
Letting $d=a^{-1/2}>1$, we obtain
\[
G\!\left(\frac32\right)
=d+1+d^{-2}+d^{-3}-3(1-d^{-2})\ln d
>d+1-3\ln d.
\]
Since $e^x\geq 1+x+x^2/2$ for $x\geq0$, substituting
$x=\ln d>0$ yields
\[
d+1-3\ln d
\geq 2-2\ln d+\frac12(\ln d)^2
=\frac12(\ln d-2)^2\geq0.
\]
Hence $G(t)\geq G(3/2)>0$ for $t\geq3/2$, so $g'(t)>0$
on $[3/2,+\infty)$.
	
 Since
\[
\frac{g(t)}{v}
=t\frac{u^{t-1}+v^{t-1}}{u^t+v^t}
\]
and $\tfrac32\leq m-\tfrac12<n-\tfrac12$, the strict
monotonicity of $g$ gives $g(m-\tfrac12)<g(n-\tfrac12)$.
Consequently,
\begin{equation}\label{E}
	\left(n-\tfrac12\right)
	\frac{u^{n-\frac32}+v^{n-\frac32}}
	{u^{n-\frac12}+v^{n-\frac12}}
	-
	\left(m-\tfrac12\right)
	\frac{u^{m-\frac32}+v^{m-\frac32}}
	{u^{m-\frac12}+v^{m-\frac12}}
	>0.
\end{equation}

On the other hand, since $E(z)=a_1(2z)^N$, we have
\[
\begin{aligned}
	&E''(u)-E''(v)
	+\bigl(E'(v)-E'(u)\bigr)
	\left(m-\tfrac12\right)
	\frac{u^{m-\frac32}+v^{m-\frac32}}
	{u^{m-\frac12}+v^{m-\frac12}}\\
	&\quad=\bigl(E'(v)-E'(u)\bigr)
	\left(
	(N-1)\frac{u^{N-2}-v^{N-2}}{v^{N-1}-u^{N-1}}
	+\left(m-\tfrac12\right)
	\frac{u^{m-\frac32}+v^{m-\frac32}}
	{u^{m-\frac12}+v^{m-\frac12}}
	\right).
\end{aligned}
\]
Since $0<u<v$,  we have the following inequalities
\[
-\frac1v<
\frac{u^{N-2}-v^{N-2}}{v^{N-1}-u^{N-1}}\leq0,
\qquad
\frac{u^{m-\frac32}+v^{m-\frac32}}
{u^{m-\frac12}+v^{m-\frac12}}>\frac1v.
\]
Thus,
$$
(N-1)\frac{u^{N-2}-v^{N-2}}{v^{N-1}-u^{N-1}}
+\left(m-\tfrac12\right)
\frac{u^{m-\frac32}+v^{m-\frac32}}
{u^{m-\frac12}+v^{m-\frac12}}>
\frac{m-N+\frac{1}{2}}{v}>0.$$
 Since $E'(v)-E'(u)=2Na_1((2v)^{N-1}-(2u)^{N-1})>0$, this inequality
together with \eqref{E} yields $F_1''(u)-F_2''(v)>0$.
\end{proof}
We now prove that $\Delta'(r_*)>0$ at every positive zero
$r_*$ of $\Delta$. Consider $r_*>0$ such that $\Delta(r_*)=0$.
By Lemma~\ref{lem2}, for each $i=1,2$, the half-length $r_*$
uniquely determines the peak height $s_i=S_i(0,r_*)$.
Let $z_i(y;s_i)$ be the corresponding solution arc of
\eqref{3}. 
{For each $i\in\{1,2\}$ and $0\leq h<s_i$, the arc
$z_i(y;s_i)$ intersects the line $z=h$ at
$(\alpha_i(h,s_i),h)$ and $(\beta_i(h,s_i),h)$.}
By 
Theorem~\ref{thm1}, the zero $r_*$ determines a unique
periodic orbit of system~\eqref{system}. Let $(0,\alpha)$
and $(0,\beta)$ be its intersections with the negative and
positive $y$-axes, respectively, where $\alpha<0<\beta$.
Thus 
\[
\alpha_1(0,s_1)=\alpha_2(0,s_2)=\alpha,
\qquad
\beta_1(0,s_1)=\beta_2(0,s_2)=\beta.
\]

Let $h_*$ be the unique positive zero of $H(z)$, given by
\[
h_*=\frac12(-a_2)^{1/(n-m)}.
\]
\medskip
\noindent\textbf{Step 1. Relationship of the relative midpoints at a common height.}
We claim that
\begin{enumerate}
	\item[\bf (a)] For all $h\geq0$ and $r>0$,
	\begin{equation}\label{12}
		\mathcal M_{F_1}(h,r)+\mathcal M_{F_2}(h,r)>0.
	\end{equation}
	
	\item[\bf (b)] For all $h\geq h_*$ and $r>0$,
	\begin{equation}\label{13}
		\mathcal M_{F_1}(h,r)>\mathcal M_{F_2}(h,r).
	\end{equation}
\end{enumerate}
Since $F_1(z)+F_2(z)=2E(z)$ and $E(z)$ is strictly increasing on
$[0,+\infty)$, for every $h\geq0$ and $z>0$ we  can obtain that
\[
\bigl(F_1(h+z)-F_1(h)\bigr)
+\bigl(F_2(h+z)-F_2(h)\bigr)
=2\bigl(E(h+z)-E(h)\bigr)>0.
\]
 It follows from Lemma~\ref{lem9}(ii)  with $h_1=h_2=h$ that
\eqref{12} holds.

Since 
\[
H(z)=(2z)^{m+\frac12}\bigl(a_2+(2z)^{n-m}\bigr),
\]
it is not difficult to obtain that 
$H(z)$ is strictly decreasing on $[0,h_0)$ and strictly increasing on
$[h_0,+\infty)$,
where \[
h_0=\frac12\left(\frac{(-a_2)(2m+1)}{2n+1}\right)^{1/(n-m)}.
\]
Hence $H(z)$ is strictly increasing
on $[h_*,+\infty)$ {because} $h_0<h_*$.  Thus, for every $h\geq h_*$
and $z>0$ we have
\[
\bigl(F_1(h+z)-F_1(h)\bigr)
-\bigl(F_2(h+z)-F_2(h)\bigr)
=2\bigl(H(h+z)-H(h)\bigr)>0.
\]
 By Lemma~\ref{lem9}(i) with $h_1=h_2=h${, we obtain} \eqref{13}.
 
 \medskip
 \noindent\textbf{{Step 2. 
    Comparison of half-widths at equal heights.}}
 { We claim that
 	$r_2(h,s_2)>r_1(h,s_1)$	for every $h\in(0,\min\{s_1,s_2\})$.}
 	As shown in the proof of
 	Lemma~\ref{lem2}, these functions are $C^1$ in $h$ for
 	fixed $s_i$. Differentiating the identities
 	\[
 	z_i(\alpha_i(h{,s_i});s_i)=h,\qquad z_i(\beta_i(h{,s_i});s_i)=h
 	\]
 	with respect to $h$ and using \eqref{3}, we obtain
 	\begin{equation}\label{ab}
 	\frac{d \alpha_i(h,s_i)}{d h}=\frac{1}{F_i(h)-\alpha_i(h,s_i)},\qquad
 		\frac{d \beta_i(h,s_i)}{d h}=-\frac{1}{\beta_i(h,s_i)-F_i(h)}.
 	\end{equation}
 
 Define
 \[
 R_1(h)=\frac{\alpha_2(h,s_2)-\alpha_1(h,s_1)}{2},
 \qquad
 R_2(h)=\frac{\beta_2(h,s_2)-\beta_1(h,s_1)}{2}.
 \]
 By the definition of the half-length \eqref{r},
 \[
 \begin{aligned}
 	r_2(h,s_2)-r_1(h,s_1)
 	&=\frac{\beta_2(h,s_2)-\alpha_2(h,s_2)}{2}
 	-\frac{\beta_1(h,s_1)-\alpha_1(h,s_1)}{2}\\
 	&=R_2(h)-R_1(h).
 \end{aligned}
 \]
 Thus, according to \eqref{ab} and $F_1(h)-F_2(h)=2H(h)$,
 we have
 \[
 \begin{aligned}
 	R_1'(h)
 	=\frac{H(h)}{K_1(h)}+\frac{R_1(h)}{K_1(h)},\qquad 
 	R_2'(h)
 	=\frac{H(h)}{K_2(h)}+\frac{R_2(h)}{K_2(h)}
 \end{aligned}
 \]
 and
 \begin{equation}\notag
 	\bigl(R_2(h)-R_1(h)\bigr)'
 =\frac{R_2(h)-R_1(h)}{K_1(h)}
 +\left(\frac{K_1(h)-K_2(h)}{K_1(h)K_2(h)}\right)
 \bigl(H(h)+R_2(h)\bigr),
 \end{equation}
 where
 \[
 \begin{aligned}
 	K_1(h)
 	&=\bigl(F_1(h)-\alpha_1(h,s_1)\bigr)
 	\bigl(F_2(h)-\alpha_2(h,s_2)\bigr)>0,\\
 	K_2(h)
 	&=\bigl(\beta_1(h,s_1)-F_1(h)\bigr)
 	\bigl(\beta_2(h,s_2)-F_2(h)\bigr)>0.
 \end{aligned}
 \]
 By the variation-of-constants formula and $R_1(0)=R_2(0)=0$,
we obtain
\[
\begin{aligned}
	R_1(h)
	=\int_0^h \frac{H(z)}{K_1(z)}
	\exp\!\left(\int_z^h \frac{dw}{K_1(w)}\right)\,dz,\qquad
	R_2(h)
=\int_0^h \frac{H(z)}{K_2(z)}
	\exp\!\left(\int_z^h \frac{dw}{K_2(w)}\right)\,dz
\end{aligned}
\]
and
\[
R_2(h)-R_1(h)
=\int_0^h
\exp\!\left(\int_z^h \frac{dw}{K_1(w)}\right)
\left(\frac{K_1(z)-K_2(z)}{K_1(z)K_2(z)}\right)
\bigl(R_2(z)+H(z)\bigr)\,dz.
\]
Since $H(z)<0$ for $z\in(0,h_*)$,  the integral formulas imply $R_1(h)<0$ and $R_2(h)<0$ whenever $0<h<\min\{s_1,s_2\}$ and $h\leq h_*$. On the other hand, by \eqref{12},
\[
\alpha+\beta
=\mathcal M_{F_1}(0,r_*)+\mathcal M_{F_2}(0,r_*)>0.
\]
Thus $K_2(0)=\beta^2>\alpha^2=K_1(0)$. By continuity,
$K_1(z)-K_2(z)<0$ for sufficiently small $z>0$. 
Together with  $R_2(z)+H(z)<0$ for such $z$, $R_2(h)-R_1(h)>0$ for sufficiently
small $h>0$.
We claim that this inequality {holds} throughout
$(0,\min\{s_1,s_2\})$. Otherwise, there is a first height
$\bar h$ on $(0,\min\{s_1,s_2\})$ such that $R_2(\bar h)=R_1(\bar h)$.
At $\bar h$, set
\[
\bar r=r_1(\bar h,s_1)=r_2(\bar h,s_2).
\]

Furthermore, by \eqref{r} and \eqref{M}, at $h=\bar h$ we have
\[
\begin{aligned}
	K_1(\bar h)-K_2(\bar h)
	&=\bigl(\bar r-\mathcal M_{F_1}(\bar h,\bar r)\bigr)
	\bigl(\bar r-\mathcal M_{F_2}(\bar h,\bar r)\bigr)-\bigl(\bar r+\mathcal M_{F_1}(\bar h,\bar r)\bigr)
	\bigl(\bar r+\mathcal M_{F_2}(\bar h,\bar r)\bigr)\\
	&=-2\bar r\bigl(\mathcal M_{F_1}(\bar h,\bar r)
	+\mathcal M_{F_2}(\bar h,\bar r)\bigr).
\end{aligned}
\]
Moreover, it follows 
from \eqref{12} that $\mathcal M_{F_1}(\bar h,\bar r)
+\mathcal M_{F_2}(\bar h,\bar r)>0$. Thus,	$K_1(\bar h)-K_2(\bar h)<0$.

If $\bar h\leq h_*$, then $R_2(\bar h)<0$ and
$H(\bar h)\leq0$. If $\bar h>h_*$, then   by \eqref{r} and \eqref{M},
$$R_2(\bar h)+H(\bar h)=\frac{\beta_2(\bar h,s_2)-\beta_1(\bar h,s_1)}{2}+\frac{F_1(\bar h)-F_2(\bar h)}{2}=\frac{
\mathcal M_{F_2}(\bar h,\bar r)+\bar r}{2}-\frac{\mathcal M_{F_1}(\bar h,\bar r)+\bar r}{2}.
$$
Moreover by \eqref{13}, it follows that $\mathcal M_{F_1}(\bar h,\bar r)>\mathcal M_{F_2}(\bar h,\bar r)$.
Thus, $R_2(\bar h)+H(\bar h)<0$ {in both cases.}
Consequently,
$$
	\bigl(R_2(\bar h)-R_1(\bar h)\bigr)'=\left(\frac{K_1(\bar h)-K_2(\bar h)}{K_1(\bar h)K_2(\bar h)}\right)
	\bigl(H(\bar h)+R_2(\bar h)\bigr)>0.
$$
However, $R_2(h)-R_1(h)>0$ for $0<h<\bar h$ and vanishes
at $\bar h$. Taking the derivative from the left therefore
gives $(R_2-R_1)'(\bar h)\leq0$, a contradiction. 
This proves that
\[
r_2(h,s_2)-r_1(h,s_1)=R_2(h)-R_1(h)>0
\]
for every $h\in(0,\min\{s_1,s_2\})$.

\medskip
\noindent\textbf{Step 3. Comparison of relative midpoints at equal
	half-lengths.}
By \eqref{ab},
$$
\frac{dr_i(h,s_i)}{dh}=-\left(\frac{1}{2(F_i(h)-\alpha_i(h,s_i))}+\frac{1}{2(\beta_i(h,s_i)-F_i(h))}\right)<0,
$$
and so $r_i(h,s_i)$ is strictly decreasing in $h$. 
Since $r_i(0,s_i)=r_*$ and $r_i(h,s_i)\to0$ as $h\to s_i^-$,
for each $r\in(0,r_*]$, there is a unique height $h_i(r)\in[0,s_i)$ such
that $r_i(h_i(r),s_i)=r$. In particular,
$h_1(r_*)=h_2(r_*)=0$.

We first show that $s_1\leq s_2$ and  $h_1(r)<h_2(r)$ for $0<r<r_*$. Indeed, if $s_2<s_1$, 
letting $h\to s_2^-$ in $r_2(h,s_2)>r_1(h,s_1)$ would
give $0\geq r_1(s_2,s_1)>0$. Thus, $s_1\leq s_2$.  By Step~2,
$$
r_2(h_1(r),s_2)>r_1(h_1(r),s_1)=r=r_2(h_2(r),s_2),
$$
so the strict decrease of $r_2(h,s_2)$ yields 
 $0<h_1(r)<h_2(r)$ for $0<r<r_*$.

Differentiating $r_i(h_i(r),s_i)=r$ with respect to $r$ {and using} \eqref{r}, \eqref{M} and \eqref{ab}, we obtain that
\begin{equation}\notag
\frac{d h_i(r)}{dr}=
-\frac{r^2-\mathcal M_{F_i}^2(h_i(r),r)}{r},
\end{equation}
 and 
\begin{equation}\label{15}
\frac{d}{dr}\mathcal M_{F_i}(h_i(r),r)
=-\frac{\mathcal M_{F_i}(h_i(r),r)}{r}
+\frac{r^2-\mathcal M^2_{F_i}(h_i(r),r)}{r}
F_i'(h_i(r)).
\end{equation}

We next claim that for all $r\in(0,r_*)$
\begin{equation}\label{14}
	\mathcal M_{F_1}(h_1(r),r)>\mathcal M_{F_2}(h_2(r),r).
\end{equation}
Define
\[
D(r)=\mathcal M_{F_1}(h_1(r),r)-\mathcal M_{F_2}(h_2(r),r).
\]
Since $h_1(r_*)=h_2(r_*)=0$ and $\Delta(r_*)=0$, we have
$D(r_*)=0$.  In the proof of this claim, we write $\mathcal M_{F_i}(r)$
for $\mathcal M_{F_i}(h_i(r),r)$, for $i=1,2$.
By \eqref{15}, we have that 
\[
\begin{aligned}
	D'(r)=
	&-\frac{1+\bigl(\mathcal M_{F_1}(r)+\mathcal M_{F_2}(r)\bigr)
		F_2'(h_2(r))}{r}\,D(r)\\
	&\qquad +\frac{r^2-\mathcal M_{F_1}^2(r)}{r}
	\bigl(F_1'(h_1(r))-F_2'(h_2(r))\bigr).
\end{aligned}
\]
Applying the variation-of-constants formula with $D(r_*)=0$,
we obtain
\begin{equation*}
\begin{aligned}
	D(r)
	&=\int_r^{r_*}
	\exp\!\left(
	\int_r^\rho
	\frac{1+\bigl(\mathcal M_{F_1}(w)+\mathcal M_{F_2}(w)\bigr)
		F_2'(h_2(w))}{w}\,dw
	\right)\\
	&\qquad{}\times
	\frac{\rho^2-\mathcal M_{F_1}^2(\rho)}{\rho}
	\bigl(F_2'(h_2(\rho))-F_1'(h_1(\rho))\bigr)\,d\rho.
\end{aligned}
\end{equation*}

For $\rho<r_*$ sufficiently close to $r_*$, we have
$0<h_1(\rho)<h_2(\rho)\ll h_0$. Since $E'$ is strictly increasing
and $H'<0$ on $(0,h_0)$,
\[
F_2'(h_2(\rho))-F_1'(h_1(\rho))
=E'(h_2(\rho))-E'(h_1(\rho))-H'(h_2(\rho))-H'(h_1(\rho))>0.
\]
Moreover, by \eqref{8}, $|\mathcal M_{F_1}(h_1(\rho),\rho)|=|\mathcal M_{F_1}(\rho)|<\rho$.
 Hence $D(r)>0$
for such $r<r_*$ sufficiently close to $r_*$.

Now suppose that $D(r)$ has a zero in $(0,r_*)$, and let
$\widehat r$ be its largest one in this interval. 
{Since $D(r)>0$ for $\hat r<r<r_*$ and $D(\hat r)=0$, we have}
$$
D'(\hat r)\geq 0.
$$
If we can show that
$F_1'(h_1(\widehat r))<F_2'(h_2(\widehat r))$, 
{then evaluating the differential equation for $D$ at $r=\hat r$, where $D(\hat r)=0$, gives}
\[
D'(\widehat r)
=\frac{\widehat r^{\,2}-\mathcal M_{F_1}^2(\widehat r)}
{\widehat r}
\bigl(F_1'(h_1(\widehat r))-F_2'(h_2(\widehat r))\bigr)<0,
\]
which yields a contradiction. Therefore, our goal is to prove that
$F_1'(h_1(\widehat r))<F_2'(h_2(\widehat r))$.
Otherwise, $F_1'(h_1(\widehat r))\geq F_2'(h_2(\widehat r))$.
Consider the function 
\[
L(z)=F_1(h_1(\widehat r)+z)-F_2(h_2(\widehat r)+z)
-\bigl(F_1(h_1(\widehat r))-F_2(h_2(\widehat r))\bigr)
\]
with $z\geq -h_1(\widehat r)$.
One can check that $L(0)=0$ and
$L'(0)=F_1'(h_1(\widehat r))-F_2'(h_2(\widehat r))\geq0$.
If $L'(z)>0$ for all $z>0$, then $L(z)>L(0)\geq 0$. 
By Lemma~\ref{lem9}(i) with $h_i=h_i(\widehat r)$, we have
\[
\mathcal M_{F_1}(h_1(\widehat r),\widehat r)
>
\mathcal M_{F_2}(h_2(\widehat r),\widehat r),
\]
which contradicts $D(\widehat r)=0$. 
Thus,  we have  $L'(z)\leq0$ for some $z>0$.
If $L'(0)>0$, continuity gives $L'(z)>0$ for sufficiently
small $z>0$. If $L'(0)=0$, by Lemma~\ref{lem10}, $L''(0)>0$. Hence  $L'(z)>0$ for sufficiently
small $z>0$   also holds. 
Consequently,
      $L'(z)$ has a
smallest positive zero $z_0$, and $L'(z)>0$ on $(0,z_0)$.
Since $0<h_1(\widehat r)+z_0<h_2(\widehat r)+z_0$, it follows from
Lemma~\ref{lem10} that $L''(z_0)>0$. On the other hand,
$L'(z)>0$ on $(0,z_0)$ and $L'(z_0)=0$ imply
$L''(z_0)\leq0$, a contradiction.  

In summary,  $D(r)$ has  no  zeros in $(0,r_*)$.
Together with $D(r)>0$
for such $r<r_*$ sufficiently close to $r_*$, it follows {that} $D(r)>0$ for the whole interval and \eqref{14} holds.

\medskip
\noindent\textbf{Step 4. The sign of $\Delta'(r_*)$.}

By the definition of $\Delta$,
\[
\Delta'(r_*)
=\left.\left(\frac{\partial }{\partial r}(\mathcal M_{F_1}(0,r))
-\frac{\partial }{\partial r}(\mathcal M_{F_2}(0,r))\right)\right|_{r=r_*}
\]
 Further, from \eqref{22}, 
$$
\left.\frac{\partial }{\partial r}(\mathcal M_{F_i}(0,r))\right|_{r=r_*}=
 \left.\frac{
 	\dfrac{\partial}{\partial s}(\beta_i(0,s))
 	+\dfrac{\partial}{\partial s}(\alpha_i(0,s))
 }{
 	\dfrac{\partial}{\partial s}(\beta_i(0,s))
 -\dfrac{\partial}{\partial s}(\alpha_i(0,s))
 }\right|_{s=s_i}
 $$
 Let
 \[
 J_i=\int_\alpha^\beta F_i'\bigl(z_i(y;s_i)\bigr)\,dy.
 \]
 Since $F_i(0)=0$, $\alpha_i(0,s_i)=\alpha$, and
 $\beta_i(0,s_i)=\beta$, it follows from \eqref{4} that
 \[
 -\left.\frac{\partial}{\partial s}(\alpha_i(0,s))\right|_{s=s_i}
 =\frac{\xi_i(\alpha;s_i)}{-\alpha},
 \qquad
 \left.\frac{\partial}{\partial s}(\beta_i(0,s))\right|_{s=s_i}
 =\frac{\xi_i(\beta;s_i)}{\beta}.
 \]
Using \eqref{xi}, we obtain
 \[
 \left.
 \frac{\dfrac{\partial}{\partial s}\beta_i(0,s)}
 {-\dfrac{\partial}{\partial s}\alpha_i(0,s)}
 \right|_{s=s_i}
 =\frac{-\alpha}{\beta}\times
 \frac{\xi_i(\beta;s_i)}{\xi_i(\alpha;s_i)}
 =\frac{-\alpha}{\beta}e^{J_i}.
 \]
Consequently,
\[
\left.\frac{\partial }{\partial r}(\mathcal M_{F_i}(0,r))\right|_{r=r_*}=\frac{-\alpha e^{J_i}-\beta}{\beta-\alpha e^{J_i}},
\]
and hence
\[
\Delta'(r_*)
=\frac{-2\alpha\beta\bigl(e^{J_1}-e^{J_2}\bigr)}
{(\beta-\alpha e^{J_1})(\beta-\alpha e^{J_2})}.
\]
Since $\alpha<0<\beta$, the sign of $\Delta'(r_*)$ is
the sign of $J_1-J_2$.

The solution $z_i(y;s_i)$ increases from $0$ to $s_i$ when 
$y$ increases from $\alpha$ to $F_i(s_i)$, and $z_i(y;s_i)$ decreases from
$s_i$ to $0$ when $y$ {increases} from $F_i(s_i)$ to $\beta$.
On the increasing and decreasing branches, the substitution
$h=z_i(y;s_i)$ gives $y=\alpha_i(h,s_i)$ and
$y=\beta_i(h,s_i)$, respectively. Thus,
\[
\begin{aligned}
	J_i
	&=\int_\alpha^{F_i(s_i)} F_i'\bigl(z_i(y;s_i)\bigr)\,dy
	+\int_{F_i(s_i)}^\beta F_i'\bigl(z_i(y;s_i)\bigr)\,dy\\
	&=\int_0^{s_i} F_i'(h)
	\frac{\partial}{\partial h}(\alpha_i(h,s_i))\,dh
	+\int_{s_i}^0 F_i'(h)
	\frac{\partial}{\partial h}(\beta_i(h,s_i))\,dh\\
	&=-2\int_0^{s_i}F_i'(h)
	\frac{\partial}{\partial h}(r_i(h,s_i))\,dh.
\end{aligned}
\]
Moreover, considering  $r=r_i(h,s_i)$, we have  $h=h_i(r)$ and
$dr=\bigl(\frac{\partial}{\partial h}r_i(h,s_i)\bigr)\,dh$.
When  $h$ increases from $0$ to $s_i$, $r$ decreases from
$r_*$ to $0$. Therefore,
\[
J_i=-2\int_{r_*}^0 F_i'\bigl(h_i(r)\bigr)\,dr
=2\int_0^{r_*}F_i'\bigl(h_i(r)\bigr)\,dr.
\]

For each fixed $i=1,2$, define
\[
T_i(r)=
\begin{cases}
	\displaystyle \frac12\ln\Big(
	\frac{r+\mathcal M_{F_i}(h_i(r),r)}
	{r-\mathcal M_{F_i}(h_i(r),r)}
	\Big), & 0<r\leq r_*,\\[6pt]
	0, & r=0.
\end{cases}
\]
We claim that $T_i\in C^1([0,r_*])$.
First, the definition of $T_i$ gives $T_i\in C^0((0,r_*])$.
 For $0<r<r_*$, \eqref{15} gives
\begin{equation}\label{T}
    T_i'(r)
    =
    \frac{
        r{\frac{d}{dr}}\bigl(M_{F_i}(h_i(r),r)\bigr)
        -M_{F_i}(h_i(r),r)
    }{
        r^2-M_{F_i}^2(h_i(r),r)
    }
    =
    F_i'(h_i(r))
    -\frac{2M_{F_i}(h_i(r),r)}
    {r^2-M_{F_i}^2(h_i(r),r)}.
\end{equation}
{Since $|M_{F_i}(0,r_*)|<r_*$ and $h_i(r_*)=0$,
$T_i'(r)$ has a finite limit as $r\to r_*^-$.
Thus $T_i\in C^1((0,r_*])$.}
Now, we focus on $r=0$.
Since $|\mathcal M_{F_i}(h_i(r),r)|<r$ and $F_i'$ is
bounded on $[0,s_i]$, we obtain from  \eqref{15}  
\[
\frac{d}{dr}\bigl[r\mathcal M_{F_i}(h_i(r),r)\bigr]
=\bigl(r^2-\mathcal M_{F_i}^2(h_i(r),r)\bigr)
F_i'(h_i(r))
=\mathcal O(r^2),
\]
so $\mathcal M_{F_i}(h_i(r),r)=\mathcal O(r^2)$.
In particular, $\mathcal M_{F_i}(h_i(r),r)/r\to0$ when $r\to 0^+$, so
$T_i(r)\to T_i(0)=0$. Therefore, $T_i(r)$ is continuous at $r=0$.  

 Since $h_i(r)\to s_i$ as $r\to0^+$, l'Hôpital's rule gives
 \[
 \begin{aligned}
 	\lim_{r\to0^+}\frac{\mathcal M_{F_i}(h_i(r),r)}{r^2}
 	&=\lim_{r\to0^+}
 	\frac{\dfrac{d}{dr}\bigl(r\mathcal M_{F_i}(h_i(r),r)\bigr)}
 	{3r^2}\\
 	&=\lim_{r\to0^+}\frac13
 	\left(1-\frac{\mathcal M_{F_i}^2(h_i(r),r)}{r^2}\right)
 	F_i'(h_i(r))
 	=\frac13F_i'(s_i).
 \end{aligned}
 \]
 Using $\frac12\ln\bigl((1+x)/(1-x)\bigr)=x+\mathcal O(x^3)$ as
 $x\to0$, with $x=\mathcal M_{F_i}(h_i(r),r)/r$, we obtain
 \[
 \lim_{r\to0^+}\frac{T_i(r)}{r}
 =\lim_{r\to0^+} \left(\frac{\mathcal M_{F_i}(h_i(r),r)}{r^2}+\mathcal O(r^2)\right)
=\frac13F_i'(s_i).
 \]
{Since $T_i(0)=0$, the preceding limit gives
\[
    T_i'(0)
    =\lim_{r\to0^+}\frac{T_i(r)-T_i(0)}{r}
    =\frac13F_i'(s_i),
\]
where the derivative at $r=0$ is understood as the right derivative.
On the other hand, by \eqref{T},
\[
    T_i'(r)
    =F_i'(h_i(r))
    -\frac{2M_{F_i}(h_i(r),r)/r^2}
    {1-M_{F_i}^2(h_i(r),r)/r^2}.
\]
Using $h_i(r)\to s_i$, $M_{F_i}(h_i(r),r)=O(r^2)$,
and
\[
    \lim_{r\to0^+}\frac{M_{F_i}(h_i(r),r)}{r^2}
    =\frac13F_i'(s_i),
\]
we obtain
\[
    \lim_{r\to0^+}T_i'(r)
    =F_i'(s_i)-\frac23F_i'(s_i)
    =\frac13F_i'(s_i)
    =T_i'(0).
\]
Hence $T_i\in C^1([0,r_*])$, and \eqref{T}
extends continuously to $[0,r_*]$.}

Finally,
by \eqref{T} and the expressions for $J_i$,
we obtain
\[
\begin{aligned}
	J_1-J_2
	={}&2\int_0^{r_*}\bigl(T_1'(r)-T_2'(r)\bigr)\,dr\\
	&+4\int_0^{r_*}\left(
	\frac{\mathcal M_{F_1}(h_1(r),r)}
	{r^2-\mathcal M_{F_1}^2(h_1(r),r)}
	-\frac{\mathcal M_{F_2}(h_2(r),r)}
	{r^2-\mathcal M_{F_2}^2(h_2(r),r)}
	\right)\,dr.
\end{aligned}
\]
The first integral vanishes because $T_1(0)=T_2(0)=0$ and
$T_1(r_*)=T_2(r_*)=\frac12\ln\bigl(\beta/(-\alpha)\bigr)$.
For each fixed $r\in(0,r_*)$, the function
$f(x)= x/(r^2-x^2)$ is strictly increasing on $(-r,r)$.
Since $|\mathcal M_{F_i}(h_i(r),r)|<r$ and $\mathcal M_{F_1}(h_1(r),r)>\mathcal M_{F_2}(h_2(r),r)$  in \eqref{14}, 
$$	\frac{\mathcal M_{F_1}(h_1(r),r)}
{r^2-\mathcal M_{F_1}^2(h_1(r),r)}
-\frac{\mathcal M_{F_2}(h_2(r),r)}
{r^2-\mathcal M_{F_2}^2(h_2(r),r)}>0.$$
 Therefore $J_1-J_2>0$,
and hence $\Delta'(r_*)>0$.

Once we have shown that $\Delta'(r_*)>0$ at every positive
zero $r_*$ of $\Delta$, it follows that $\Delta$ has at most
one positive zero, since $\Delta$ is $C^1$. By
Theorem~\ref{thm1}, the corresponding limit cycle is
hyperbolic and asymptotically stable. The proof of Theorem \ref{thm2} is finished.

	\section*{Acknowledgements}
	The work was supported by  Hunan Basic Science Research Center for Mathematical Analysis (Grant No. 2024JC2002) and by
	the  National Natural Science Foundation of China (Grant Nos. 12322109 and 12571190). 	 The fourth author is supported by the Postdoctoral Fellowship Program and China Postdoctoral Science Foundation under Grant {No.} BX2026012.
	
	{\bf Data Availability:}  {No datasets were generated or analysed during this study.}

{\bf Declaration:} The authors have no conflicts of interest to declare that are relevant to the content of this article.


\end{document}